\documentclass[10pt]{article}
\usepackage{amssymb}
\usepackage{latexsym}
\input amssym.def
\input amssym
\newtheorem{theorem}{Theorem}[section]
\newtheorem{lemma}[theorem]{Lemma}

\newtheorem{corollary}[theorem]{Corollary}

\newtheorem{remark}[theorem]{Remark}

\newcommand{\ds}{\displaystyle}

\newcommand{\halmos}{\rule{1ex}{1.4ex}}

\newcommand{\bea}{\begin{eqnarray}}
\newcommand{\epf}{\hspace*{\fill}\mbox{$\halmos$}}
\newcommand{\eea}{\end{eqnarray}}
\newcommand{\nn}{\nonumber \\}
\newcommand{\be}{\begin{equation}}

\newcommand{\ee}{\end{equation}}

\usepackage{amssymb}

\title{{Ramanujan's ``Lost Notebook'' and the Virasoro Algebra} }
\author{Antun  Milas}
\begin{document}
\date{}
\bibliographystyle{alpha}
\maketitle
\begin{abstract}
\noindent By using the theory of vertex operator algebras, we gave
a new proof of the famous Ramanujan's modulus 5 modular equation
from his "Lost Notebook" (p.139 in \cite{R}). Furthermore, we
obtained an infinite list of $q$-identities for all odd moduli;
thus, we generalized the result of Ramanujan.
\end{abstract}

\renewcommand{\theequation}{\thesection.\arabic{equation}}
\renewcommand{\thetheorem}{\thesection.\arabic{theorem}}
\setcounter{equation}{0}
\setcounter{theorem}{0}
\setcounter{section}{0}

\section{Introduction}
\noindent According to Hardy (cf. \cite{A}, p.177): ``..If I had
to select one formula from all Ramanujan's work, I would agree
agree with Major MacMahon in selecting...
\be \label{mod5} \sum_{n \geq 0} p(5n+4)q^n=\frac{5
(q^5;q^5)_\infty^5}{(q;q)_\infty^6}, \ee where $p(n)$ is the
number of partitions of $n$'', and
$$(a;q)_\infty=(1-a)(1-aq) \cdots.$$
Closely related to formula (\ref{mod5}) is a pair of
$q$--identities recorded by Ramanujan in his ``Lost Notebook''
(cf. p.139-140, \cite{R}): \be \label{rr1} 1-5 \sum_{n \geq 0}
\left(\frac{n}{5}\right) \frac{n q^n}{1-q^n}=
\frac{(q;q)_\infty^5}{(q^5;q^5)_\infty} \ee and \be \label{rr2}
\sum_{n \geq 1} \left( \frac{n}{5} \right) \frac{q^n}{(1-q^n)^2}
=\frac{q (q^5;q^5)^5_{\infty}}{(q;q)_\infty}, \ee where
$\left(\frac{n}{5} \right)$ is the Legendre symbol.
As a matter of fact, we can talk about a single identity (cf. \cite{C2}) because
(\ref{rr2}) (and subsequently (\ref{mod5}))
can be obtained easily from (\ref{rr1}) by applying a classical result of Hecke
(cf. p. 119 in \cite{My}).
By now, there are several proofs of (\ref{rr1}) and (\ref{rr2}) in
the literature. The first proof was given by Bailey (\cite{B1} and
\cite{B2}) by using the ${}_6 \Psi_6$--summation. For recent
proofs see, for instance, \cite{Rg}, \cite{C1} and references
therein. For an extensive account on Ramanujan's modular
identities see \cite{Br} (see also \cite{A1}, \cite{BrO}).

Compared to these conventional approaches 
(e.g.,  hypergeometric $q$--series, modular forms), our approach to 
Ramanujan's modular identities is based on completely 
different ideas. Let us elaborate on recent
developments and results that brought us to Ramanujan's ``Lost
Notebook'' \cite{R} and in particular to (\ref{mod5}).

It is well known that infinite-dimensional Lie theoretical methods
can be used to conjecture, interpret and ultimately prove series
of combinatorial and $q$--series identities related to partitions.
This direction was initiated by Lepowsky and Wilson in \cite{LW} and it is based
on explicit constructions of integrable highest weight representations for affine Lie algebras.
In addition, various dilogarithm techniques, crystal bases and path representation techniques that originate in conformal field theory and statistical physics, led to interesting combinatorial and 
$q$--series identities. 
%
Besides affine Lie algebras there is another important algebraic structure closely related
to affine Lie algebras: the Virasoro algebra (cf. \cite{FFu},
\cite{FFu1}, \cite{KR}, \cite{KW}). Even though the Virasoro algebra and its representations
theory are well-understood (including
character formulas and their modular properties), the 
``smallest'' representations of the Virasoro
algebra (i.e, the minimal models \cite{BPZ}) have no known
explicit constructions.
%

Our motivation for studying Ramanujan's ``Lost Notebook''
identities stems from the following fact: A large part of
Ramanujan's work concerns modular equations closely related to
some of his $q$--identities and continued fractions. Up to now
there has not been any work done in the direction of understanding
these identities from conformal field theoretical point of view.
This is surprising because the {modular invariance} of
characters holds for a large class of vertex operator
algebras (VOA) \cite{Zh}. The Virasoro algebra is already included
in the definition of VOA, so it appears very natural to seek
for modular identities in connections with irreducible Virasoro
algebra modules (e.g., minimal models).

Let us briefly outline the content of the paper.
In the first part we consider the simplest (yet quite involved)
minimal models with exactly two irreducible modules and with 
$c=\frac{-22}{5}$ (i.e., the Lee-Yang model). This model
is related to Rogers-Ramanujan identities. The main idea is to show that the
irreducible characters satisfy a second order linear differential
equation with coefficients being certain Eisenstein series. In
order to achieve this we use the theory of vertex
operator algebras (especially  Zhu's work \cite{Zh}). When we
combine the character formulas for the Virasoro minimal models
obtained in \cite{FFu}--\cite{FFu2} with some standard ODE
techniques to obtain (\ref{rr1}) (cf. Theorem \ref{mod5lost}).

In the second part we provide a generalization of the formula
(\ref{rr1}). As in the $c=\frac{-22}{5}$ case, the key idea is to
consider a series of ODEs satisfied by irreducible characters. It
is a highly nontrivial problem to compute these ODEs explicitly. Luckily, for our present purposes, it was enough to
compute only the first two leading coefficients. As a consequence, we obtain the following family of $q$-identities (the $k=2$ case
corresponds to Ramanujan's modular equation (\ref{rr1})):
\smallskip
\\
\noindent Fix $k \in \mathbb{N}$, $k \geq 2$. For every $i=1,...,k$, let
$$A_{i}(q)=\frac{6i^2-6i+1+6k^2-12ki+5k}{12(2k+1)}+
\sum_{{n \geq 0,}{ n \neq \pm i,0 \ {\rm mod} \ (2k+1)}} \frac{nq^n}{1-q^n}.$$
Then
\bea \label{detintro}
&& \left| \begin{array}{ccccc} 1 & 1 & . & . & 1 \\
                              A_1(q) & A_{2}(q) & . & . & A_{k}(q) \\
                              \bar{P}_2(A_1(q)) & \bar{P}_2(A_2(q)) & . & . &  \bar{P}_2(A_k(q)) \\
                              . & . & . & . &  . \\
                              \bar{P}_{k-1}(A_1(q)) & \bar{P}_{k-1}(A_2(q)) & . & . & \bar{P}_{k-1}(A_k(q))
\end{array} \right| \nn
&& =\frac{\ds{\prod_{i=1}^{k-1}} (2i)!}{(-4k-2)^{\frac{k(k-1)}{2}
} }
\left(\frac{(q;q)^{2k+1}_\infty}{(q^{2k+1};q^{2k+1})_\infty}\right)^{k-1},
\eea where $\bar{P}_j$'s are certain (shifted) Fa\`a di Bruno
operators (polynomials) defined in (\ref{faashift}). We should mention that
our identities do not have obvious modular properties (at least
not the determinant side), so it is an open question to express
(\ref{detintro}) in a more explicit way (for related work see
\cite{M2}).

\noindent It appears that our modulus 7 identity ($k=3$, see Corollary \ref{moduli7}) might be related to another pair of modular identities recorded by
Ramanujan (cf. p.145 in \cite{R}). We shall treat a possible connection
in our future publications.

{\bf Acknowledgments:} It was indeed hard to trace all the known
proofs of (\ref{mod5}), (\ref{rr1}) and (\ref{rr2}). We apologize if some
important references are omitted. We would like to thank Jim
Lepowsky for conversations on many related subjects. A few years
ago Lepowsky and the author were trying to relate classical
Rogers-Ramanujan identities and Zhu's work \cite{Zh}. We also
thank Bruce Berndt for pointing us to \cite{BrO} and
Steve Milne for bringing \cite{Mi} to our attention.

\renewcommand{\theequation}{\thesection.\arabic{equation}}
\renewcommand{\thetheorem}{\thesection.\arabic{theorem}}
\setcounter{equation}{0}
\setcounter{theorem}{0}

\section{The Virasoro algebra and minimal models}

The Virasoro algebra ${\rm Vir}$ (cf. \cite{FFu}, \cite{KR}, etc.) is defined
as the unique non-trivial central extension of the Lie
algebra of polynomial vector fields on $\mathbb{C}^*$. It is
generated by $L_n$, $n \in \mathbb{Z}$ and $C$, with bracket relations \be
[L_m,L_n]=(m-n)L_{m+n}+\frac{m^3-m}{12}\delta_{m+n,0} C, \ee where
$\delta_{m+n,0}$ is the Kronecker symbol and $C$ is the central
element. Let us fix a triangular decomposition
$${\rm Vir}={\rm Vir}_+ \oplus {\rm Vir}_0 \oplus {\rm Vir}_-,$$
where
${\rm Vir}_+$ is spanned by $L_i$, $i > 1$, ${\rm Vir}_-$ is spanned
by $L_i$, $i<0$ and ${\rm Vir}_0$ is spanned by $C$ and $L_0$.
Let $M$ be a ${\rm Vir}$--module. We shall denote the action of $L_n$ on
$M$ by $L(n)$. Let $U({\rm Vir})$ denote the enveloping algebra of
${\rm Vir}$. Thus
$$U({\rm Vir})=U({\rm Vir}_-) \otimes U({\rm Vir}_0) \otimes U({\rm Vir}_+).$$
The enveloping algebra $U({\rm Vir}_-)$ is equipped with the natural
filtration
$$\mathbb{C}=U({\rm Vir}_-)_0 \subset U({\rm Vir}_-)_1 \subset \cdots.$$
It follows from PBW theorem that
every element of $U({\rm Vir}_-)_k$ is spanned by
elements of the form
$$L(-i_1) \cdots L(-i_r), \ \ 0 \leq r \leq k, \ i_j >0, \ j=1,...,r.$$
We shall denote by ${\rm Vir}_{\leq -n}$ (resp. ${\rm Vir}_{\geq n}$)
a Lie subalgebra spanned by $L_i$, $i \leq -n$ (resp. $i \geq n$).

Let $c,h \in \mathbb{C}$. Let $\mathbb{C} v_{c,h}$ denote a one-dimensional
$U( {\rm Vir}_{\geq 0})$--module such that
$$L(n) v_{c,h}=0, \ \ n>0,$$
$$C \cdot v_{c,h}=c v_{c,h}, \ \ L(0) v_{c,h}=h v_{c,h}.$$
Consider the Verma module \cite{KR} \be M(c,h)=U({\rm Vir})
\otimes_{U({\rm Vir}_+)} \mathbb{C}v_{c,h}, \ee We shall say that
$M(c,h)$ has {\em central charge} $c$ and {\em weight} $h$. In
particular $M(c,h)$ has the maximal submodule $M^{(1)}(c,h)$ and
the corresponding irreducible quotient
$$L(c,h)=M(c,h)/M^{(1)}(c,h).$$
In 1980s Feigin and Fuchs  provided a detailed embedding structure
for Verma modules for all values of $c$ and $h$ \cite{FFu},
\cite{FFu1}.

There is an infinite, distinguished, family of irreducible highest weight
modules $L(c_{p.q},h^{m,n}_{p,q})$ ({\em minimal models})
parameterized by the central charge
$$c_{p,q}=1-\frac{6(p-q)^2}{pq},$$
where $p, q \in \mathbb{N}$, $ p,q \geq 2$, $(p,q)=1$, and
with weights
$$h^{m,n}_{p,q}=\frac{(np-mq)^2-(p-q)^2}{4 pq},$$
where $1 \leq m <p$, $1 \leq n <q$. Notice that for certain
pairs $(m,n)$ and $(m',n')$, 
$$h^{m,n}_{p,q}=h^{m',n'}_{p,q}.$$
More precisely, there are exactly
$$\frac{(p-1)(q-1)}{2}$$
different values of $h_{m,n}$ for $1 \leq m <p$, $1 \leq n <q$.
Because of
$$M(c,h) \cong U({\rm Vir}_-)$$
and the fact that $M^{(1)}(c,h)$ is graded it is clear that
$L(c_{p,q},h_{m,n})$ is naturally
$\mathbb{Q}$-graded with respect to the action
of $L(0)$. Moreover, the graded subspaces are finite--dimensional.
Hence to every highest weight module $M$ we can associate its
graded dimension, $q$--trace, or simply its {\em character}
\be
{\rm tr}|_M q^{L(0)},
\ee
where (unless otherwise stated) $q$ is just a formal variable \footnote{Formal variable $q$ has nothing to do with the integer $q$ used for parameterization of the central charge.}.
In the case of minimal models we shall write
$$ {\rm ch}_{c_{p.q},h^{m,n}_{p,q}}(q)={\rm tr}|_{L(c_{p,q},h^{m.n}_{p,q})} q^{L(0)}.$$
It is not hard to see that \be {\rm tr}|_{M(c,h)}
q^{L(0)}=\frac{q^h}{(q;q)_\infty}. \ee However, computing ${\rm
tr}|_{L(c,h)} q^{L(0)}$ is a much more difficult problem. The only
known proof uses a complete BGG-type resolution for irreducible
highest weight modules due to Feigin and Fuchs \cite{FFu}. By
using their result it is a straightforward task to obtain explicit
formulas for \be \label{character} {\rm
ch}_{c_{p.q},h^{m,n}_{p,q}}(q), \ee where $p,q,m$ and $n$ are as
above.
We should mention that some partial result have been known
prior to their result. 

For present purposes the expression (\ref{character}) has to be
modified (actually, this modification turns out to be essential).
Let \be \label{modchar} \bar{{\rm
ch}}_{c_{p.q},h^{m,n}_{p,q}}(q)={\rm
tr}|_{L(c_{p.q},h^{m,n}_{p,q})} q^{L(0)-\frac{c_{p,q}}{24}}. \ee
From now on we will consider a one-parameter family of central
charges
$$c_{2,2k+1}=1-\frac{6(2k-1)^2}{(4k+2)}, \ \ k \geq 2$$
and the corresponding weights
$$h^{1,i}_{2,2k+1}=\frac{(2(k-i)+1))^2-(2k-1)^2}{8(2k+1)}, \ \ i=1,...,k.$$
The $c=c_{2,3}$ case is not interesting because it gives the trivial module.

An important observation (cf. \cite{RC}) is that
characters of minimal models with the central charge $c_{2,2k+1}$, $k \geq 1$, can be
expressed as infinite products  (see \cite{FFr}, \cite{RC}, \cite{KW}).
We will explore this fact in the later sections.

\renewcommand{\theequation}{\thesection.\arabic{equation}}
\renewcommand{\thetheorem}{\thesection.\arabic{theorem}}
\setcounter{equation}{0}
\setcounter{theorem}{0}

\section{The $c=c_{2,5}$ case}
The simplest non-trivial minimal models occur for
$c=c_{2,5}=-\frac{22}{5}$. In this case there are exactly two
irreducible modules:
$$L\left(\frac{-22}{5},0\right) \ \ \mbox{and} \ \  L\left(\frac{-22}{5},\frac{-1}{5}\right).$$
The corresponding
characters (written as infinite products) are essentially product sides
appearing in Rogers-Ramanujan identities \cite{FFr}, \cite{RC}.
More precisely (cf. \cite{KW}, \cite{FFr}, \cite{RC}):
\be \label{first} {\rm ch}_{-22/5,0}(q)=\prod_{n \geq 0}
\frac{1}{(1-q^{5n+2})(1-q^{5n+3})} \ee and
\be \label{second}
{\rm ch}_{-22/5,-1/5}(q)=q^{-1/5} \prod_{n \geq 0} \frac{1}
{(1-q^{5n+1})(1-q^{5n+4})}. \ee
We will show that these character formulas can be used
to obtain Ramanujan's modular equation mentioned in the introduction.

\renewcommand{\theequation}{\thesection.\arabic{equation}}
\renewcommand{\thetheorem}{\thesection.\arabic{theorem}}
\setcounter{equation}{0}
\setcounter{theorem}{0}

\section{Vertex operator algebras and modular invariance}
In this section we shall recall some of results from the theory of vertex operator algebras.


For the definition of vertex operator algebras, modules for vertex operator algebras
and irreducible modules see \cite{FHL} or \cite{FLM}. It is well known (cf.
\cite{FZ}, \cite{W}) that the so--called {\em vacuum}
module
$$V(c,0)=M(c,0)/\langle L(-1){v_{c,0}} \rangle$$
can be equipped with a vertex operator algebra
structure such that
$${\bf 1}=v_{c,0}$$
and
$$\omega=L(-2){\bf 1}.$$
By quotienting $V(c,0)$ by the maximal ideal we obtain a
simple vertex operator algebra $L(c,0)$. However, $L(c,0)$ is
not very interesting for all values of $c$. In the 
$c=c_{2,2k+1}$ case (and more generally $c=c_{p,q}$), the representation
theory of $L(c_{2,2k+1},0)$ becomes surprisingly simple (cf. \cite{W}, see also \cite{FZ}).
\begin{theorem} \label{classmod}
For every $k \geq 1$, the vertex operator algebra $L(c_{2,2k+1},0)$ is rational
(in the sense of \cite{DLM1}). Moreover, a complete
list of (inequivalent) irreducible $L(c_{2,2k+1},0)$-modules is
given by
$$L(c_{2,2k+1},h^{1,i}_{2,2k+1}), \ \ i=1,...,k.$$
In particular, the only irreducible
$L\left(\frac{-22}{5},0\right)$--modules are (up to isomorphism)
$L\left(\frac{-22}{5},0\right)$ and
$L\left(\frac{-22}{5},\frac{-1}{5}\right)$.
\end{theorem}
The previous result is a reformulation in the
language of vertex operator algebras of a result due to Feigin and
Fuchs \cite{FFu2}.

\renewcommand{\theequation}{\thesection.\arabic{equation}}
\renewcommand{\thetheorem}{\thesection.\arabic{theorem}}
\setcounter{equation}{0}
\setcounter{theorem}{0}

\section{A change of ``coordinate'' for vertex operator algebras}
Let $V$ be an arbitrary vertex operator algebra and suppose that
$u \in V$ is a homogeneous element (i.e., an eigenvector for $L(0)$). Let
$$Y[u,y]=e^{y {\rm deg}(u)} Y(u,e^{y}-1),$$
where
$$L(0) u=({\rm deg}(u)) u,$$
$y$ is a formal variable and $(e^y-1)^{-n-1}$, $n \in \mathbb{Z}$, is expanded inside
$\mathbb{C}((y))$, truncated Laurent series in $y$.
Extend $Y[u,y]$ to all $u \in V$ by the linearity. Let
$$Y[u,y]=\sum_{n \in \mathbb{Z}} u[n]y^{-n-1}, \ u[n] \in {\rm End}(V).$$
The following theorem was proven in \cite{H} (see also \cite{Zh} and \cite{L1}).
\begin{theorem} \label{isomorphism}
Let $(V,Y( \cdot,y), {\omega},{\bf 1})$ be a vertex operator
algebra and
$$\tilde{\omega}=L[-2]{\bf 1}=(L(-2)-\frac{c}{24}){\bf 1} \in V.$$
The quadruple $(V,Y[ \cdot,y], \tilde{\omega},{\bf 1})$ has a
vertex operator algebra structure isomorphic to the vertex operator
algebra $(V,Y( \cdot,y), {\omega},{\bf 1}).$
In particular, if we let
$$Y[L[-2]{\bf 1},x]=\sum_{n \in \mathbb{Z}} L[n]x^{-n-2}$$
then
$$[L[m],L[n]]=(m-n)L[m+n]+\frac{m^3-m}{12}\delta_{m+n,0}c.$$
\end{theorem}
The following lemma is from \cite{Zh} (see also Chapter 7 in \cite{H}).
\begin{lemma}
For every $n \geq 0$ there are sequences
$$\{ c_n^{(i)} \}, \ c_n^{(i)} \in \mathbb{Q}, \ i \geq (n+1)$$
and
$$\{ d_n^{(i)} \}, \ d_n^{(i)} \in \mathbb{Q}, \ i \geq (n+1),$$
such that
\be
L(n)=L[n]+ \sum_{i \geq (n+1)} c_{n}^{(i)} L[i]
\ee
and
\be
L[n]=L(n)+ \sum_{i \geq (n+1)} d_{n}^{(i)} L(i).
\ee
\end{lemma}
Now, we specialize $V=V(c,0)$.
The following lemma is essentially from \cite{Zh}.
\begin{lemma} \label{singform}
Let $v \in V(c,0)$ be a singular vector, i.e., a homogeneous vector
annihilated by $L(i)$, $i>0$.
Suppose that
$$v=\sum a_{I} L(-i_1)L(-i_2) \cdots L(-i_k){\bf 1},$$
where $a_{I} \in \mathbb{C}$, and the summation goes over all indices
$i_1 \geq i_2 \geq \cdots \geq i_k $, such that $i_1+ \cdots +i_k={\rm deg}(v).$
Then
$$v=\sum a_{I} L[-i_1]L[-i_2] \cdots L[-i_k]{\bf 1}.$$
Informally speaking, singular vectors in the vertex operator algebra
$V(c,0)$ are invariant with respect to the change of coordinate $x \mapsto e^x-1$.
\end{lemma}
{\em Proof:} The proof is a consequence of Huang's theorem
\cite{H} concerning an arbitrary change of coordinate induced
by a conformal transformation vanishing at $0$ (see \cite{H} for
details, cf. also \cite{Zh}). Let $\Psi_{e^x-1}$ be the isomorphism
of $V$ induced by the change of variable $x \longrightarrow
e^x-1$. Let us compute $\Psi_{e^x-1}(v)$.
$$v=v_{-1}{\bf 1} \mapsto e^{\sum_{i \geq 1} r_i L_i} v=
\Psi_{e^x-1} v=v[-1]{\bf 1}.$$
for some $r_i \in \mathbb{C}$.
Every singular vector $v$ satisfies
$$L(i) v=0, \ \ i \geq 1.$$
Hence
$$v \mapsto v, \ ({\rm under} \ \ \Psi_{e^x-1}).$$
Because of the previous corollary (in fact Theorem
\ref{isomorphism}) $v$ is also a singular vector with respect to
$L[i]$ generators. On the other hand, by the definition of
isomorphism for VOA,
$$\Psi_{e^x-1} : V \longrightarrow V, \ \ \Psi_{e^x-1}(Y(u,x)v)=Y[\Psi_{e^x-1}(u),x] \Psi_{e^x-1}(v).$$
Therefore
$$\Psi_{e^x-1} (L(-i_1) \cdots L(-i_k){\bf 1})=L[-i_1] \cdots L[-i_k]{\bf 1},$$
for any choice  of $i_1,...,i_k$.
\epf

\begin{remark}
{\em It is important to mention  that the previous construction
has been known by physicists since early eighties (after the
seminal work \cite{BPZ}). Invariance of singular vectors (or {\em
primary fields}) with respect to conformal transformations is one
of the most important features in {conformal field theory}.}
\end{remark}

\noindent Because of $L[-2]{\bf 1}=(L(-2)-\frac{c}{24}){\bf 1}$, it is
convenient to introduce
$$\bar{L}(0)=L(0)-\frac{c}{24}.$$
This transformation corresponds to cylindrical change of coordinates.
Notice also that $L[0] \neq L(0)-\frac{c}{24}$. The following theorem
is essentially due to Zhu \cite{Zh} (for further generalizations
and modifications see \cite{DLM}).
\begin{theorem}
Let $V$ be a rational vertex operator algebra which satisfies the
$C_2$--condition \footnote{$C_2$ condition: The vector space
spanned by $u_{-2}v$,
$u,v \in V$, has a finite codimension \cite{DLM}.}. Let
$M_1$,...$M_k$ be a list of all (inequivalent) irreducible
$V$--modules. Then the vector space spanned by
$${\rm tr}|_{M_1} q^{\bar{L}(0)},...,{\rm tr}|_{M_k} q^{\bar{L}(0)},$$
is modular invariant with respect to $\Gamma(1)$,
where $\gamma$ acts on the modulus $\tau$ ($q=e^{2 \pi i \tau}$) in the
standard way
$$ \gamma \cdot \tau=\frac{a\tau+b}{c \tau+d}, \ \ \gamma=\left[\begin{array}{cc} a & b \\ c & d \end{array} \right] \in \Gamma(1).$$
\end{theorem}
Now, we let  $V=L \left(\frac{-22}{5},0 \right)$, where
$$\bar{L}(0)=L(0)+\frac{11}{60}.$$
The previous theorem implies the following result (even though we will not
use it in the rest of the paper).
\begin{corollary} \label{coro}
The vector space spanned by
\be \label{firstmod} \bar{{\rm
ch}}_{-22/5,0}(q)=q^{11/60} \prod_{n
\geq 0} \frac{1} {(1-q^{5n+2})(1-q^{5n+3})}  \ee and \be
\label{secondmod} \bar{{\rm ch}}_{-22/5,-1/5}(q)= q^{-1/60} \prod_{n \geq 0} \frac{1}
{(1-q^{5n+1})(1-q^{5n+4})} \ee
is modular invariant.
\end{corollary}
Let $\tilde{G}_{2k}(q)$, $k \geq 1$, denote (normalized)
Eisenstein series given by their $q$--expansions
$$\tilde{G}_{2k}(q)=\frac{-B_{2k}}{(2k)!}+\frac{2}{(2k-1)!}\sum_{n \geq
0} \frac{n^{2k-1} q^n}{1-q^n},$$
where $B_{2k}$, $k \geq 1$, are Bernoulli numbers. In particular,
$$\tilde{G}_{2}(q)=\frac{-1}{12}+2 \sum_{n \geq
0} \frac{n q^n}{1-q^n}$$
and
$$\tilde{G}_{4}(q)=\frac{1}{720}+\frac{1}{3} \sum_{n \geq
0} \frac{n^3 q^n}{1-q^n}.$$ Our normalization is convenient
because all the coefficients in the $q$--expansion of
$\tilde{G}_{2k}(q)$ are rational numbers (notice that Zhu \cite{Zh}
used a different normalization, cf. \cite{M1}).

Let $V$ be a vertex operator algebra and $M$ a $V$--module.
Also, let $u \in V$ be a homogeneous elements. Define
$$o(u)=u_{{\rm wt}-1} \in {\rm End}(M).$$
For instance
$$o(\omega)=L(0).$$
Extend this definition for every $u \in V$ by
linearity. Also,
$$o(\tilde{\omega})=L(0)-\frac{c}{24}.$$
The
following result was proven in \cite{Zh} (also it is a consequence
of a more general result obtained in \cite{M1}; see also \cite{DMN}).
\begin{theorem} \label{recursion}
For every $u,v \in V$,
\be \label{zero}
{\rm tr}|_M o(u[0]v)q^{\bar{L}(0)}=0
\ee
and
\bea \label{recursiony}
&& {\rm tr}|_M o(u[-1]v) q^{\bar{L}(0)} \nn
&& = {\rm tr}|_M o(u)o(v)q^{\bar{L}(0)}+ \sum_{k \geq 1}
\tilde{G}_{2k}(q){\rm tr}|_M X(u[2k-1]v,x)q^{\bar{L}(0)}.
\eea
\end{theorem}

\renewcommand{\theequation}{\thesection.\arabic{equation}}
\renewcommand{\thetheorem}{\thesection.\arabic{theorem}}
\setcounter{equation}{0}
\setcounter{theorem}{0}

\section{Differential equation}
In this section we obtain a second order linear differential equation with a fundamental system of solutions given
by $\bar{{\rm ch}}_{-22/5,0}(q)$ and $\bar{{\rm ch}}_{-22/5,-1/5}(q)$.
We should mention that Kaneko and Zagier (cf. \cite{KZ})
considered related second order
differential equations from a different point of view.
\begin{theorem} \label{2}
Let $\tilde{G}_2(q)$ and $\tilde{G}_4(q)$ be as above and $q=e^{2
\pi i \tau}$, $\tau \in \mathbb{H}$. Then $\bar{ch}_{-22/5,0}(q)$ and $\bar{ch}_{-22/5,-1/5}(q)$
form a fundamental system of solutions of
 \be \label{diffmain} \left(q
\frac{d}{dq} \right)^2 F(q)+2 \tilde{G}_2(q)\left(q \frac{d}{dq}
\right) F(q)-\frac{11}{5}\tilde{G}_4(q) F(q)=0.
\ee
\end{theorem}
{\em Proof:} It is enough to show that both
$\bar{ch}_{-22/5,0}(q)$ and $\bar{ch}_{-22/5,-1/5}(q)$ satisfy the
equation (\ref{diffmain}). Firstly, from the structure of Verma
modules for the Virasoro algebra \cite{FFu} it follows that
$$v=(L^2(-2)-\frac{3}{5} L(-4)){\bf 1}$$
is a singular vectors inside $V(-22/5,0)$.
This vector generates the maximal submodule
of the vertex operator algebra $V(-22/5,0)$, i.e.,
$$L(\frac{-22}{5},0)=V(\frac{-22}{5},0)/ \langle v \rangle,$$
where $\langle S \rangle$ denotes the ${\rm Vir}$--submodule generated by
the set $S$. By using Lemma 
\ref{singform}, it follows that
$$v=(L^2[-2]-\frac{3}{5} L[-4]){\bf 1}.$$
Let $M$ be an arbitrary $L(-22/5,0)$--module. Then
\be \label{null} Y_M((L^2[-2]-\frac{3}{5} L[-4]){\bf 1},x)=0 \ee
inside
$${\rm End}(M)[[x,x^{-1}]].$$
Hence
\be \label{use}
{\rm tr}|_M o((L^2[-2]-\frac{3}{5} L[-4]){\bf 1})q^{L(0)}=0.
\ee
In particular, we may set
$M=L\left(\frac{-22}{5},0 \right)$ or
$M=L\left(\frac{-22}{5},\frac{-1}{5} \right)$.
Now, we apply the formula (\ref{zero}) and get
$${\rm tr}|_M o(\tilde{\omega}[0]v)q^{\bar{L}(0)}={\rm tr}|_M o(L[-1]v) q^{L(0)}=0,$$
for every $v \in V$. We shall pick $v=L[-3]{\bf 1}$, which implies
$${\rm tr}|_M o(L[-1]L[-3]{\bf 1}) q^{L(0)}=2 {\rm tr}|_M o(L[-4]{\bf 1}) q^{L(0)}=0.$$
The previous formula and (\ref{use}) imply \be \label{difftwo} {\rm tr}|_M
o(L^2[-2]{\bf 1},x) q^{L(0)}=0. \ee
From
$${\rm tr}|_M o(L[-2]L[-2]{\bf 1})q^{\bar{L}(0)}=
{\rm tr}|_M o(\tilde{\omega}[-1]L[-2]{\bf 1})q^{\bar{L}(0)}$$
and
$$o(L[-2]{\bf 1})=L(0)+\frac{11}{60}$$
we get
\bea
&& {\rm tr}|_M o(L[-2]L[-2]{\bf 1})q^{\bar{L}(0)} \nn
&& = {\rm tr}|_M o(L[-2]{\bf 1})o(L[-2]{\bf 1}) q^{\bar{L}(0)}+ 2 \tilde{G}_2(q)
{\rm tr}|_M o(L[-2]{\bf 1})q^{\bar{L}(0)} \nn
&& -\frac{11}{5} \tilde{G}_4(q) {\rm tr}|_M q^{\bar{L}(0)} \nn
&& ={\rm tr}|_M (L(0)-\frac{c_{2,5}}{24})^2 q^{\bar{L}(0)}
+2\tilde{G}_2(q)
{\rm tr}|_M (L(0)-\frac{c_{2,5}}{24})q^{\bar{L}(0)} \nn
&& -\frac{11}{5} \tilde{G}_4(q) {\rm tr}|_M q^{\bar{L}(0)} \nn
&&= \left(q \frac{d}{dq} \right)^2 {\rm tr}|_M q^{\bar{L}(0)}
+2 \tilde{G}_2(q) \left(q \frac{d}{dq} \right) {\rm tr}|_M q^{\bar{L}(0)}-
\frac{11}{5}\tilde{G}_4(q){\rm tr}|_M q^{\bar{L}(0)}.
\eea
\epf
\begin{remark}
{\em The property (\ref{difftwo}) is closely related to
"difference--two condition at the distance one" property (cf.
\cite{A}). More precisely, Feigin and Frenkel \cite{FFr} used (\ref{null})
(and the sewing rules obtained in \cite{BFM}) to give a conformal
field theoretical proof of Rogers--Ramanujan identities and their
generalizations (due to Gordon and Andrews).}
\end{remark}
Let
$$\left( \frac{n}{5} \right) $$
be the Legendre symbol. In his Lost Notebook (p. 139, \cite{R}), S. Ramanujan
recorded
\begin{theorem} \label{mod5lost}
\be \label{lost}
1-5 \sum_{n \geq 0} \left( \frac{n}{5} \right) \frac{n q^n}{1-q^n}
=\frac{ (q;q)^5_\infty }{(q^5,q^5)_\infty}.
\ee
\end{theorem}
We will need a simple lemma concerning infinite products
and their logarithmic derivatives.
\begin{lemma} \label{lemmaprod}
Let
$$\mathcal{A}(q)=q^r \prod_{n \geq 0} \frac{1}{(1-q^n)^{a_n}},$$
 where
$a_n \in \mathbb{Z}$ and $r \in \mathbb{C}$. Then
$$\left(q\frac{d}{dq}\right) \mathcal{A}(q)=
\left(r+\sum_{n \geq 0} a_n \frac{nq^n}{1-q^n} \right)\mathcal{A}(q)$$
or (equivalently)
$$\left(\frac{d}{d\tau}\right) \mathcal{A}(\tau)=
2 \pi i \left(r+\sum_{n \geq 0} a_n \frac{nq^n}{1-q^n} \right)\mathcal{A}(\tau).$$
\end{lemma}
{\em Proof of the Theorem:} The proof will follow from the following elementary (but
fundamental) result due to Abel \cite{Hi}. Let $\Omega \subset \mathbb{C}$ be a domain
and $P_1(z)$ and $P_2(z)$ be two holomorphic functions inside $D$. Suppose that
$y_1$ and $y_2$ form a fundamental system of solutions of the differential
equation
$$y''+P_1(z)y'+P_2(z)y=0.$$
Then we have the formula
\be \label{abel}
W(z)=W(z_0) e^{\ds{-\int_{z_0}^z P_1(t)dt}},
\ee
where
$$W(z)=\left| \begin{array}{cc} y_1 & y_2 \\ y'_1 & y'_2 \end{array}\right|=y_1y'_2-y'_1y_2$$
is the Wronskian of the system, $z_0 \in \Omega$, and the
integration goes along any rectifiable path in $\Omega$. We will
apply Abel's formula (\ref{abel}) to (\ref{diffmain}), i.e.,
\be \label{diffmain1}
\left(q
\frac{d}{dq} \right)^2 F(q)+2 \tilde{G}_2(q)\left(q \frac{d}{dq}
\right) F(q)-\frac{11}{5}\tilde{G}_4(q) F(q)=0.
\ee
Firstly, notice that our differential equation (\ref{diffmain1}) can be written in terms of $\tau$, rather
than $q$, where we can take $\Omega$ to be the upper-half plane and
$$'=\frac{1}{2 \pi i} \frac{d}{d \tau}=\left(q\frac{d}{dq}\right).$$
Now, for the
fundamental system of solutions of (\ref{diffmain1}) we pick (cf. Theorem \ref{2})
$$y_1(\tau)=\bar{{\rm ch}}_{-22/5,0}(\tau) \ \ {\rm and} \ \
y_2(\tau)=\bar{{\rm ch}}_{-22/5,-1/5}(\tau)
$$ It is easy to compute
the Wronskian by using the infinite product expressions for
$\bar{{\rm ch}}_{-22/5,0}(\tau)$ and
$\bar{{\rm ch}}_{-22/5,-1/5}(\tau)$ (see formulas (\ref{firstmod}) and (\ref{secondmod})).
We have
\bea \label{wronski} &&
y_1(\tau)y'_2(\tau)-y'_1(\tau)y_2(\tau)= \nn && =
\biggl\{\frac{-1}{5}+\frac{11}{60}+\sum_{n \geq 0} \left(
\frac{(5n+1)q^{5n+1}}{1-q^{5n+1}}
+\frac{(5n+4)q^{5n+4}}{1-q^{5n+4}} \right)-\frac{11}{60} \nn &&
-\sum_{n \geq 0} \left( \frac{(5n+2)q^{5n+2}}{1-q^{5n+2}}
+\frac{(5n+3)q^{5n+3}}{1-q^{5n+3}} \right) \biggr\}
y_1(\tau)y_2(\tau). \eea
By
combining (\ref{wronski}) and (\ref{diffmain1}) together with the
formula (\ref{abel}) we get \bea \label{00} &&
\biggl\{\frac{-1}{5}+\sum_{n \geq 0} \biggl(
\frac{(5n+1)q^{5n+1}}{1-q^{5n+1}}
+\frac{(5n+4)q^{5n+4}}{1-q^{5n+4}} \nn &&
-\frac{(5n+2)q^{5n+2}}{1-q^{5n+2}}
-\frac{(5n+3)q^{5n+3}}{1-q^{5n+3}} \biggr) \biggr\}\frac{\eta(5
\tau)}{\eta(\tau)}
  \nn && = W(\tau_0) e^{ \ds{-2
\int_{\tau_0}^{\tau} \tilde{G}_2(\tau) d(2 \pi i  \tau)} },
\eea
where we used the fact
$$H_1(\tau)H_2(\tau)=\frac{q^{1/6}(q^5;q^5)_\infty}{(q;q)_\infty}=\frac{\eta(5 \tau)}{\eta(\tau)}.$$
By Lemma \ref{lemmaprod}
\be
\label{eta}
\left(\frac{1}{2 \pi i} \frac{d}{d \tau}\right)
\eta^4(\tau)=
\left(\frac{4}{24}-4 \sum_{n \geq 1} \frac{nq^n}{1-q^n}\right)
\eta^4(\tau)= -2 \tilde{G}_2(\tau) \eta^4(\tau)
\ee
or
$$\eta^4(\tau)=e^{\ds{-2 \int_{\tau_0}^\tau \tilde{G}_2(t) d (2 \pi i \tau)}}.$$
The previous formula implies \be \label{22} W(\tau_0) e^{ \ds{-2
\int_{\tau_0}^{\tau} \tilde{G}_2(\tau) d (2 \pi i \tau)} }= W(\tau_0)
e^{\ds{\int^{\tau}_{\tau_0} (\frac{1}{6} - 4 \sum_{n \geq
0} \frac{n q^n}{1-q^n}) d (2 \pi i \tau) }}=C \eta^4(\tau), \ee where $C$ is
some constant which does not depend on $\tau$. Now, (\ref{00})
and (\ref{22}) imply that \bea && \sum_{n \geq 0} \biggl(\frac{-1}{5}+
\frac{(5n+1)q^{5n+1}}{1-q^{5n+1}}
+\frac{(5n+4)q^{5n+4}}{1-q^{5n+4}}
-\frac{(5n+2)q^{5n+2}}{1-q^{5n+2}}
-\frac{(5n+3)q^{5n+3}}{1-q^{5n+3}} \biggr)\frac{\eta(5
\tau)}{\eta(\tau)} \nn && = {C} \eta^4(\tau). \eea
%
Therefore \bea \label{etaq} && \frac{-1}{5}+ \sum_{n \geq 0}
\left(\frac{(5n+1)q^{5n+1}}{1-q^{5n+1}}
+\frac{(5n+4)q^{5n+4}}{1-q^{5n+4}}
-\frac{(5n+2)q^{5n+2}}{1-q^{5n+2}}
-\frac{(5n+3)q^{5n+3}}{1-q^{5n+3}} \right)= \nn && = \frac{-1}{5}
+\sum_{n \geq 0} \left(\frac{n}{5}\right)\frac{n q^n}{1-q^n}={C}
\frac{(q;q)^5_\infty}{(q^5;q^5)_\infty}. \eea By comparing the
first terms on both sides we get
$$C=\frac{-1}{5}.$$
If we multiply (\ref{etaq}) by $-5$ we get \be 1-5 \sum_{n \geq 0}
\left( \frac{n}{5} \right) \frac{n q^n}{1-q^n} =\frac{
(q;q)^5_\infty }{(q^5;q^5)_\infty}. \ee \epf

\renewcommand{\theequation}{\thesection.\arabic{equation}}
\renewcommand{\thetheorem}{\thesection.\arabic{theorem}}
\setcounter{equation}{0}
\setcounter{theorem}{0}

\section{A recursion formula}
In this section, as a byproduct of Theorem \ref{mod5lost}, we obtain two
recursion formulas for coefficients in the $q$--expansions of the Rogers-Ramanujan's $q$--series \footnote{These formulas are inefficient for
computation though.}.

Let $b(n)$ denote the number of partitions  of $n$ in parts of
the form $5i+1$ and $5i+4$, $i \geq 0$, and let $a(n)$ denote the number
of partitions of $n$ in parts of the form $5i+2$ and $5i+3$, $i
\geq 0$. If we recall Rogers-Ramanujan identities
\cite{A}, $b(n)$ is the number of partitions of $n$ satisfying the
``difference two condition at the distance one'' and $a(n)$ is the
number of partitions of $n$ satisfying the ``difference two condition at the distance one'' with the smallest part $> 1$.
Let
$$\sigma_k(n)=\sum_{d | n} d^{k}.$$
Clearly (cf. \cite{A})
$$L_1(q)=\sum_{n \geq 0} a(n) q^n =
\prod_{n \geq 0} \frac{1}{(1-q^{5n+2})(1-q^{5n+3})}=\sum_{n \geq 0} \frac{q^{n^2+n}}{(q)_n}$$
and
$$L_2(q)=\sum_{n \geq 0} b(n) q^n=
\prod_{n \geq 0} \frac{1}{(1-q^{5n+1})(1-q^{5n+4})}=\sum_{n \geq 0} \frac{q^{n^2}}{(q)_n}.$$
\begin{theorem}
For every $n \geq 1$,
\be \label{aa}
a(n)=\frac{11 \ds{\sum_{k \geq 0}^{n-1}}
\left(\sigma_3(n-k)-\sigma_1(n-k)\right)
a(k)-60 \ds{ \sum_{k \geq 0}^{n-1} } k a(k) \sigma_1(n-k)}{15n^2+3n},
\ee
\be \label{bb}
b(n)=\frac{\ds{ \sum_{k \geq 0}^{n-1} } \left(11 \sigma_3(n-k)+\sigma_1(n-k)\right)b(k)-60\ds{\sum_{k \geq 0}^{n-1}} k b(k) \sigma_1(n-k)}{15n^2-3n}.
\ee
\end{theorem}
{\em Proof:}  From the differential equation
(\ref{diffmain}) and
$$\bar{{\rm ch}}|_{-22/5,0}(q)=q^{11/60}L_1(q), \ \ \bar{{\rm ch}}|_{-22/5,-1/5}(q)=
q^{-1/60}L_2(q)$$
we obtain a pair of second order differential equations satisfied by
$L_1(q)$ and $L_2(q)$. These differential equations are explicitly given by:
\bea \label{l1}
&& \left(q \frac{d}{dq} \right)^2 F(q)+\left(\frac{1}{5}+4 \sum_{n \geq 1} \sigma_1(n) q^n
\right)\left(q \frac{d}{dq} \right)F(q)+ \nn
&& \frac{11}{15} \left(\sum_{n \geq 1} \left( \sigma_1(n) q^n- \sigma_3(n)q^n \right) \right)F(q)=0,
\eea
with a solution being $L_1(q)$, and
\bea \label{l4}
&& \left(q \frac{d}{dq} \right)^2 F(q)+\left(\frac{-1}{5}+4 \sum_{n \geq 1} \sigma_1(n) q^n
\right)\left(q \frac{d}{dq} \right)F(q) \nn
&&- \frac{1}{15} \left(\sum_{n \geq 1} \sigma_1(n) q^n + 11 \sum_{ n \geq 1}
\sigma_3(n) q^n \right)F(q)=0,
\eea
with a solution being $L_2(q)$. From these differential equations
and  initial conditions
$$L_1(0)=1, \ \ L'_1(0)=0 \ \ {\rm and} \ \ L_2(0)=1, \ \ L'_2(0)=1,$$
we can compute $a(n)$'s and $b(n)$'s by taking ${\rm Coeff}_{q^n}$ in
(\ref{l1}) and (\ref{l4}), respectively. This gives
formulas (\ref{aa}) and (\ref{bb}). \epf

\renewcommand{\theequation}{\thesection.\arabic{equation}}
\renewcommand{\thetheorem}{\thesection.\arabic{theorem}}
\setcounter{equation}{0}
\setcounter{theorem}{0}

\section{The general case}
In this part we generalize Ramanujan's modulus 5 identity
for all odd moduli.

Our starting point are certain infinite 
products that appear in Gordon--Andrews' identities \cite{A} 
(a generalization of Rogers-Ramanujan identities for odd moduli).
These infinite products are given by
$$\prod_{n  \neq \pm i,0 \ {\rm
mod}(2k+1)} \frac{1}{(1-q^{n})},$$ where $i=1,...,k$.
It is known (cf. \cite{FFu}--\cite{FFu1}, \cite{FFr}, \cite{KW})
that these expressions are closely related
to graded dimensions of minimal models with
the central charge $c_{2,2k+1}$, $k \geq 2$. More precisely,
$$\bar{{\rm ch}}_{c_{2,2k+1},h^{1,i}_{2,2k+1}}(q)
=q^{h^{1,i}_{2,2k+1}-\frac{c_{2,2k+1}}{24}} \prod_{n  \neq \pm i,0 \ {\rm
mod}(2k+1)} \frac{1}{(1-q^{n})},$$ where $i=1,...,k$. Let
us recall (cf. Theorem \ref{classmod}) that there are exactly $k$
(non-equivalent) irreducible modules for the vertex operator
algebra $L(c_{2,2k+1},0)$. Let us {\em multiply}
the modified characters:
\be \label{ppp}
\prod_{i=1}^k \bar{{\rm ch}}_{c_{2,2k+1},h^{1,i}_{2,2k+1}}(q)=
q^{{\sum_{i=1}^{k}(h^{1,1}_{2,2k+1}-c_{2,2k+1}/24)}}
\left(\frac{(q^{2k+1};q^{2k+1})_\infty}{(q;q)_\infty}\right)^{k-1}.
\ee
Miraculously,
\begin{lemma} \label{etak}
\be
\sum_{i=1}^k \left(h^{1,i}_{2,2k+1}-\frac{c_{2,2k+1}}{24} \right)
=\frac{k(k-1)}{12}=\frac{2k(k-1)}{24}.
\ee
\end{lemma}
The previous lemma implies
\be \label{productmod}
q^{{\sum_{i=1}^{k}(h^{1,i}_{2,2k+1}-c_{2,2k+1}/24)}}
\left(\frac{(q^{2k+1};q^{2k+1})_\infty}{(q;q)_\infty}\right)^{k-1}=
\left(\frac{\eta((2k+1)\tau)}{\eta(\tau)}\right)^{k-1}.
\ee
This expression indicates that the product of all modified
characters exhibits nice modular properties  
for every $c_{2,2k+1}$ value.
The next lemma can be obtained by straightforward computation (via Vandermonde determinant formula).
It will be useful for computation of the constant factor for our higher
moduli identities.
\begin{lemma} \label{constk}
Denote by
$$\bar{h}^{1,i}_{2,2k+1}=h^{1,i}_{2,2k+1}-\frac{c_{2,2k+1}}{24}=
\frac{6i^2-6i+1+6k^2-12ki+5k}{12(2k+1)},$$
where $i=1,...,k$.
Then
$$\left| \begin{array}{ccccc} 1 & 1 & . & . & 1 \\
                              \bar{h}^{1,1}_{2,2k+1} & \bar{h}^{1,2}_{2,2k+1} & . & . &
\bar{h}^{1,k}_{2,2k+1} \\
                              . & . & . & . &  . \\
                              . & . & . & . &  . \\
                             (\bar{h}^{1,1}_{2,2k+1})^{k-1} & (\bar{h}^{1,2}_{2,2k+1})^{k-1} & . & . & (\bar{h}^{1,k}_{2,2k+1})^{k-1}
\end{array} \right|=\frac{\ds{\prod_{i=1}^{k-1} (2i)!}}{(-4k-2)^{\frac{k(k-1)}{2} } }.
$$
\end{lemma}
In order to obtain differential equations with a fundamental system of solutions being
$\bar{{\rm ch}}_{c_{2,2k+1},h^{1,i}_{2,2k+1}}(q)$, $i=1,...,k$, we need a precise information regarding singular vectors in
the vertex operator algebra $V(c_{2,2k+1},0)$.

The following lemma will be crucial for our considerations.
Feigin and Frenkel exploited this fact in \cite{FFr} to obtain an upper
bound for the characters expressed as sum sides in Gordon-Andrews identities \cite{A}.
\begin{lemma} \label{singvector}
For every $k \geq 2$, the module $V(c_{2k+1},0)$ contains a singular vector of
degree $2k$ of the form
$$v_{sing,2k+1}=\left(L^{k}[-2]+\cdots \right) \cdot {\bf 1} ,$$
where the dots denote the lower order terms with
respect to the filtration of $U({\rm Vir}_{\leq -2})$.
\end{lemma}
{\em Proof:} It follows directly from Lemma \ref{singform}
and the description of singular vectors in \cite{FFu} (see \cite{FFr} for application in our situation).
\epf

The previous Lemma implies that
\be \label{singkkkla}
Y_M(v_{sing,2k+1},x)=Y_M(\left(L^{k}[-2]+\cdots \right) \cdot {\bf 1},x)=0,
\ee
for every $L(c_{2,2k+1},0)$--module $M$.
\begin{lemma} \label{kk1}
The condition (\ref{singkkkla}) yields a degree $k$ homogenuous linear differential equation
\be \label{diffkfor}
\left(q\frac{d}{dq}\right)^kF(q) + k(k-1) \ \tilde{G}_2(q)
\left(q \frac{d}{dq} \right)^{k-1} F(q) + \cdots + P_{k}(\tau) F(q) =0,
\ee
with a fundamental system of solutions being
$$ \bar{\rm ch}_{c_{2,2k+1},h^{1,i}_{2,2k+1}}(q), \ \ i=1,...,k.$$
\end{lemma}
{\em Proof:} The idea is similar as in the case of $c=\frac{-22}{5}$. However, this
time we are unable to obtain explicit formulas for differential equations; rather we obtained
first two leading derivatives, which is enough for our purposes.
The existence of a homogeneous differential equation of degree $k$ satisfied by
$${\rm tr}|_M q^{\bar{L}(0)},$$
$M$ being a $L(c_{2,2k+1},0)$--module, was already proven in \cite{Zh}. 
Therefore we only need to analyze the $(k-1)$-st coefficient in
\be
\left(q\frac{d}{dq} \right)^k {\rm tr}|_M q^{\bar{L}(0)}+
A_1(q) \left(q\frac{d}{dq} \right)^{k-1}{\rm tr}|_M q^{\bar{L}(0)} + \cdots =0.
\ee
An important observation is that $A_1(q)$ is independent of the lower order terms
in (\ref{singkkkla}). This  can be easily seen
by using the formula (\ref{recursiony}) or \cite{Zh} (cf. \cite{DMN}, \cite{M1}).
We will be using Theorem \ref{recursiony} repeatedly
to compute
$${\rm tr}|_M o(L^k[-2]{\bf 1})q^{\bar{L}(0)}.$$
Firstly,
\be \label{sss}
{\rm tr}|_M o(L^k[-2]{\bf 1})q^{\bar{L}(0)}=
\left(q \frac{d}{dq}\right)^{k} {\rm tr}|_M q^{\bar{L}(0)}+\cdots,
\ee
where dots involve lower order derivatives of ${\rm tr}|_M q^{\bar{L}(0)},$
multiplied (possibly) with certain Eisenstein series.
We will prove by the induction on $k$ that
$$A_1(q)=k(k-1)\tilde{G}_2(q), \ k \geq 2.$$
For $k=2$ this is true (cf. Theorem \ref{mod5lost}).
We compute
\bea \label{aaa}
&& {\rm tr}|_M o(L^k[-2]{\bf 1})q^{\bar{L}(0)} \\
&& = {\rm tr}|_M o(L[-2]{\bf 1}) o(L^{k-1}[-2]{\bf 1})
q^{\bar{L}(0)}+ 2(k-1) \tilde{G}_2(q) {\rm tr}|_M o(L^{k-1}[-2]{\bf 1})q^{\bar{L}(0)}+ \cdots ,
\nonumber
\eea
where the dots denote terms that do not contribute to $A_1(q)$.
Because of (\ref{sss}), the second term contributes to $A_1(q)$ with
$$2(k-1)\tilde{G}_2(q).$$
We shall work out the second term in (\ref{aaa})
\bea \label{ddd}
&& {\rm tr}|_M o(L[-2]{\bf 1}) o(L^{k-1}[-2]{\bf 1}) q^{\bar{L}(0)} =
{\rm tr}|_M o(L^{k-1}[-2]{\bf 1})o(L[-2]{\bf 1})q^{\bar{L}(0)} \nn
&& = \left(q \frac{d}{dq}\right) \left\{{\rm tr}|_M o(L^{k-1}[-2]{\bf 1})
q^{\bar{L}(0)} \right\}.
\eea
Now, if we use the induction hypothesis
\bea
&& {\rm tr}|_M o(L^{k-1}[-2]{\bf 1}) q^{\bar{L}(0)} \nn
&& = \left(q \frac{d}{dq}\right)^{k-1}
{\rm tr}|_M q^{\bar{L}(0)}+ (k-1)(k-2) \tilde{G}_2(q) \left(q \frac{d}{dq} \right)^{k-2}
{\rm tr}|_M q^{\bar{L}(0)}+ \cdots , \nonumber
\eea
where the dots denote terms with derivatives of ${\rm tr}|_M q^{\bar{L}(0)}$
being less or equal to
$(k-3)$.
By combining the previous equation and (\ref{ddd})
we obtain
$$A_1(q)=(2(k-1)+(k-1)(k-2))\tilde{G}_2(q)=k(k-1)\tilde{G}_2(q).$$
This proves the lemma.
\epf

\noindent For every $n \in \mathbb{N}$, define a nonlinear
differential operators $P_n(\cdot \ )$ in the following way:
$$\left(f(g(q))\right)^{[n]}=P_n(g(q))f^{[n]}(g(q)),$$
where $(h(q))^{[n]}:=\left(q \frac{ d}{dq}\right)^n h(q)$, for any
functions $f(q)$ and $g(q)$. For instance,
$$P_{1}(g(q))=\left(q\frac{d}{dq}\right) g(q)$$
and
$$P_2(g(q))=\left(\left(q\frac{d}{dq}\right) g(q) \right)^2+
\left(q\frac{d}{dq}\right)^2 g(q).$$ By using the Fa\`a di Bruno formula we get
$$P_n( \cdot \ )=\sum_{i_1,...,i_n}\frac{n!}{i_1! \cdots i_n !}
\left(\frac{1}{1!} \left(q \frac{d}{dq}\right)^1 ( \cdot \ )
\right)^{i_1} \cdots \left(\frac{1}{n!} \left(q
\frac{d}{dq}\right)^n ( \cdot \ ) \right)^{i_n},$$ where the summation
goes over all the $n$--tuples $i_1,...,i_n \geq 0$ such that
$$n=i_1+2i_2+ \cdots + n i_n.$$ We will need certain {\em shifted}
Fa\`a di Bruno operators which we define as
\be \label{faashift}
\bar{{P}}_n(\cdot \ )=\sum_{i_1,...,i_n} \frac{n!}{i_1! \cdots i_n !}
\left(\frac{1}{1!} \left( \cdot \ \right)^1 \right)^{i_1} \cdots
\left(\frac{1}{n!} \left(q \frac{d}{dq}\right)^{n-1} ( \cdot \ )
\right)^{i_n},
\ee
where, $n \geq 1$, and again the summation goes over
all the $n$--tuples $i_1,...,i_n \geq 0$, such that $n=i_1+2i_2+
\cdots + n i_n$. For instance
$$\bar{P}_1(f(q))=f(q)$$
and
$$\bar{P}_2(f(q))=(f(q))^2+\left(q \frac{d}{dq}\right)f(q).$$
\begin{lemma} \label{wrlemma}
Fix $k \geq 2$. For every $i=1,...,k$, let
$$y_i(\tau)=\bar{{\rm ch}}_{c_{2,2k+1},h^{1,i}_{2,2k+1}}(q),$$
and let
\be \label{xxx}
W(y_1,...,y_n)=\left| \begin{array}{ccccc} y_1(\tau) & y_2(\tau) & . & . & y_k(\tau) \\
                              y'_1(\tau) & y'_2(\tau) & . & . & y'_k(\tau) \\
                              y''_1(\tau) & y''_2(\tau) & . & . &  y''_k(\tau) \\
                              . & . & . & . &  . \\
                              y^{(k-1)}_1(\tau) & y^{(k-1)}_2(\tau) & . & . & y^{(k-1)}_k(\tau)
\end{array} \right|
\ee
be the Wronskian associated to $\{y_1,...,y_k\}$.
Here $$'=\frac{1}{2 \pi i} \frac{d}{d \tau}=\left(q \frac{d}{dq} \right).$$
Then
$$W(y_1,...,y_n)=\left(\prod_{i=1}^k y_i(\tau) \right)
\left| \begin{array}{ccccc} 1 & 1 & . & . & 1 \\
                              A_1(q) & A_2(q) & . & . & A_k(q) \\
                              \bar{P}_2(A_1(q)) & \bar{P}_2(A_2(q)) & . & . &  \bar{P}_2(A_k(q)) \\
                              . & . & . & . &  . \\
                              \bar{P}_{k-1}(A_1(q)) & \bar{P}_{k-1}(A_2(q)) & . & . &
                              \bar{P}_{k-1}(A_k(q))
\end{array} \right|,$$
where
$$A_i(q)=\frac{y'_i(\tau)}{y_i(\tau)}, \ i=1,...,k.$$
\end{lemma}
{\em Proof:}
We know that every $y_i(\tau)$ admits an infinite product form, which implies, because of
Lemma \ref{lemmaprod}, that
$$y_i(\tau)=e^{\ds{\int_{\tau_0}^\tau A_i(\tau) d (2 \pi i \tau)}}.$$
The Fa\`a di Bruno formula now gives
$$y^{(j)}_i(\tau)={P}_{j}\left(\int_{\tau_0}^\tau A_i(\tau)d \tau \right) e^{\ds{\int_{\tau_0}^\tau A_j(\tau)
 d (2 \pi i \tau)}}=\bar{{P}}_{j}(A_i(\tau)) y_i(\tau)$$
for every $j=1,...,k$. Finally, observe that in (\ref{xxx}) we can factor
$y_i$ from the $i$th column. This explains the multiplicative factor $\prod_{i=1}^k y_i(\tau)$
and proves the lemma.
\epf

Here is our main result:
\begin{theorem} \label{genk}
Fix $k \geq 2$. For every $i=1,...k$, let
$$A_{i}(q)=\bar{h}_{1,i}+\sum_{n \geq 0, n \neq \pm i,0 \ {\rm mod} \ (2k+1)} \frac{nq^n}{1-q^n}.$$
Then
\bea
&& \left| \begin{array}{ccccc} 1 & 1 & . & . & 1 \\
                              A_1(q) & A_{2}(q) & . & . & A_{k}(q) \\
                              \bar{P}_2(A_1(q)) & \bar{P}_2(A_2(q)) & . & . &  \bar{P}_2(A_k(q)) \\
                              . & . & . & . &  . \\
                              \bar{P}_{k-1}(A_1(q)) & \bar{P}_{k-1}(A_2(q)) & . & . & \bar{P}_{k-1}(A_k(q))
\end{array} \right| \nn
&& =\frac{\ds{\prod_{i=1}^{k-1}} (2i)!}{(-4k-2)^{\frac{k(k-1)}{2} } }
\left(\frac{(q;q)^{2k+1}_\infty}{(q^{2k+1};q^{2k+1})_\infty}\right)^{k-1}.
\eea
\end{theorem}
{\em Proof:}
As in the $k=2$ case
we apply Abel's theorem for the $k$th order linear differential equation (\ref{diffkfor}) in
Lemma \ref{kk1}. The same lemma implies that for a fundamental system
of solutions we can take
$$H_i(\tau)=\bar{ch}_{c_{2,2k+1},h^{1,i}_{2,2k+1}}(\tau), \ \ i=1,...,k.$$
With this choice
$$W(H_1(\tau),...,H_k(\tau))=C e^{-\ds{\int_{\tau_0}^\tau k(k-1) \tilde{G}_2(\tau) d(2 \pi i \tau)}},$$
where $C$ is a constant which does not depend on $\tau$. The last expression is by Lemma \ref{lemmaprod}
equal to
$$C \eta(\tau)^{2k(k-1)}.$$
After we apply Lemma \ref{wrlemma} we get
\bea
&& \left(\prod_{i=1}^k H_i(\tau) \right)
\left| \begin{array}{ccccc} 1 & 1 & . & . & 1 \\
                              A_1(q) & A_2(q) & . & . & A_k(q) \\
                              \bar{P}_2(A_1(q)) & \bar{P}_2(A_2(q)) & . & . &  \bar{P}_2(A_k(q)) \\
                              . & . & . & . &  . \\
                              \bar{P}_{k-1}(A_1(q)) & \bar{P}_{k-1}(A_2(q)) & . & . &
                              \bar{P}_{k-1}(A_k(q)) \end{array} \right| \nn
&&= C \eta(\tau)^{2k(k-1)}.
\eea
Formulas (\ref{ppp})--(\ref{productmod}) imply
\bea
&& \left(\frac{\eta((2k+1)\tau)}{\eta(\tau)}\right)^{k-1}\left| \begin{array}{ccccc} 1 & 1 & . & . & 1 \\
                              A_1(q) & A_2(q) & . & . & A_k(q) \\
                              \bar{P}_2(A_1(q)) & \bar{P}_2(A_2(q)) & . & . &  \bar{P}_2(A_k(q)) \\
                              . & . & . & . &  . \\
                              \bar{P}_{k-1}(A_1(q)) & \bar{P}_{k-1}(A_2(q)) & . & . &
                              \bar{P}_{k-1}(A_k(q)) \end{array} \right| \nn
&&= C \eta(\tau)^{2k(k-1)}. \nonumber
\eea
Hence
$$\left| \begin{array}{ccccc} 1 & 1 & . & . & 1 \\
                              A_1(q) & A_2(q) & . & . & A_k(q) \\
                              \bar{P}_2(A_1(q)) & \bar{P}_2(A_2(q)) & . & . &  \bar{P}_2(A_k(q)) \\
                              . & . & . & . &  . \\
                              \bar{P}_{k-1}(A_1(q)) & \bar{P}_{k-1}(A_2(q)) & . & . &
                              \bar{P}_{k-1}(A_k(q)) \end{array} \right|=
C \left(\frac{\eta(\tau)^{2k+1}}{\eta((2k+1)\tau)}\right)^{k-1}.$$
To figure the constant $C$ we use Lemma \ref{constk}. The proof follows.
\epf
\vskip 5mm

\renewcommand{\theequation}{\thesection.\arabic{equation}}
\renewcommand{\thetheorem}{\thesection.\arabic{theorem}}
\setcounter{equation}{0}
\setcounter{theorem}{0}

\section{Example: a modulus 7 identity}
Here, we derive
a $q$--identity in the $c_{2,7}=-\frac{68}{7}$ case.
There are three (inequivalent) minimal models:
$$L \left(\frac{-68}{7},0\right), \ \  \left(\frac{-68}{7},\frac{-2}{7}\right) \ \mbox{and} \
\left(\frac{-68}{7},\frac{-3}{7}\right).$$ If we apply Theorem \ref{genk} we get
\begin{corollary} \label{moduli7}
Let $$'=q\frac{d}{dq}$$
and
$${A_1}(q)=\frac{17}{42} + \sum_{n \geq 0, n = 2,3,4,5 \ {\rm mod} \  7}
\frac{n q^n}{1-q^n},$$
$${A_2}(q)=\frac{5}{42} + \sum_{n \geq 0, n = 1,3,4,6 \ {\rm mod} \ 7}
\frac{n q^n}{1-q^n},$$
$${A_3}(q)=\frac{-1}{42} + \sum_{n \geq 0, n = 1,2,5,6 \ {\rm mod} \ 7}
\frac{n q^n}{1-q^n}.$$
Then
\be \label{3}
\left| \begin{array}{ccc} 1 & 1 & 1 \\ {A_1}(q) & {A_2}(q) & {A_3}(q) \\
{A_1}'(q) & {A_2}'(q) & {A_3}'(q) \end{array} \right|+
\left| \begin{array}{ccc} 1 & 1 & 1 \\ {A_1}(q) & {A_2}(q) & {A_3}(q) \\
{A_1}^2(q) & {A_2}^2(q) & {A_3}^2(q) \end{array} \right|=-\frac{6}{7^3} \left(\frac{(q;q)^{7}_\infty}{(q^7;q^7)_\infty}\right)^2.
\ee
\end{corollary}
\renewcommand{\theequation}{\thesection.\arabic{equation}}
\renewcommand{\thetheorem}{\thesection.\arabic{theorem}}
\setcounter{equation}{0}
\setcounter{theorem}{0}

\section{Future work}
\noindent (a)
After we finished the first draft of the paper, S. Milne pointed to us that
his recent work \cite{Mi} might be related to our determinantal identities.
It would be nice to understand this more precisely but perhaps in the framework of
vertex operator superalgebras (e.g. for $N=1$ and  $N=2$
superconformal models). Zhu's work \cite{Zh} (cf. \cite{DLM}) indicates that $C_2$--condition
implies existence of certain differential equation so, hopefully, one can obtain many interesting modular
identities. \\
\noindent (b) (Added in the final version) The methods of this
paper can be extended to all $c_{p,q}$--series \cite{M2}. Our main
result in \cite{M2} is an extension of certain Dyson-Macdonald's
identities.

\vskip 5mm
\noindent {\em \small \sc Department of Mathematical Sciences,
Rensselaer Polytechnic Institute, Troy, NY 12180 } \\
{\em E-mail address}: milasa@rpi.edu

\end{document}